\documentclass{article}

\usepackage{amssymb,latexsym,amsmath}

\usepackage{graphicx}

\usepackage{hyperref}

\begin{document}

\newcommand{\bfi}{\bfseries\itshape}

\makeatletter

\@addtoreset{figure}{section}

\def\thefigure{\thesection.\@arabic\c@figure}

\def\fps@figure{h, t}

\@addtoreset{table}{bsection}

\def\thetable{\thesection.\@arabic\c@table}

\def\fps@table{h, t}

\@addtoreset{equation}{section}

\def\theequation{\thesubsection.\arabic{equation}}

\makeatother

\newtheorem{theorem}{Theorem}[section]

\newtheorem{proposition}[theorem]{Proposition}

\newtheorem{lema}[theorem]{Lemma}

\newtheorem{corollary}[theorem]{Corollary}

\newtheorem{definition}[theorem]{Definition}

\newtheorem{remark}[theorem]{Remark}

\newtheorem{exempl}{Example}[section]

\newenvironment{exemplu}{\begin{exempl}  \em}{\hfill $\square$

\end{exempl}}

\newcommand{\comment}[1]{\par\noindent{\raggedright\texttt{#1}

\par\marginpar{\textsc{Comment}}}}

\newcommand{\todo}[1]{\vspace{5 mm}\par \noindent \marginpar{\textsc{ToDo}}\framebox{\begin{minipage}[c]{0.95 \textwidth}

\tt #1 \end{minipage}}\vspace{5 mm}\par}

\newcommand{\ea}{\mbox{{\bf a}}}

\newcommand{\eu}{\mbox{{\bf u}}}

\newcommand{\ueu}{\underline{\eu}}

\newcommand{\ueo}{\overline{u}}

\newcommand{\oeu}{\overline{\eu}}

\newcommand{\ew}{\mbox{{\bf w}}}

\newcommand{\ef}{\mbox{{\bf f}}}

\newcommand{\eF}{\mbox{{\bf F}}}

\newcommand{\eC}{\mbox{{\bf C}}}

\newcommand{\en}{\mbox{{\bf n}}}

\newcommand{\eT}{\mbox{{\bf T}}}

\newcommand{\eL}{\mbox{{\bf L}}}

\newcommand{\eR}{\mbox{{\bf R}}}

\newcommand{\eV}{\mbox{{\bf V}}}

\newcommand{\eU}{\mbox{{\bf U}}}

\newcommand{\ev}{\mbox{{\bf v}}}

\newcommand{\eve}{\mbox{{\bf e}}}

\newcommand{\uev}{\underline{\ev}}

\newcommand{\eY}{\mbox{{\bf Y}}}

\newcommand{\eK}{\mbox{{\bf K}}}

\newcommand{\eP}{\mbox{{\bf P}}}

\newcommand{\eS}{\mbox{{\bf S}}}

\newcommand{\eJ}{\mbox{{\bf J}}}

\newcommand{\eB}{\mbox{{\bf B}}}

\newcommand{\eH}{\mbox{{\bf H}}}

\newcommand{\leb}{\mathcal{ L}^{n}}

\newcommand{\eI}{\mathcal{ I}}

\newcommand{\eE}{\mathcal{ E}}

\newcommand{\hen}{\mathcal{H}^{n-1}}

\newcommand{\eBV}{\mbox{{\bf BV}}}

\newcommand{\eA}{\mbox{{\bf A}}}

\newcommand{\eSBV}{\mbox{{\bf SBV}}}

\newcommand{\eBD}{\mbox{{\bf BD}}}

\newcommand{\eSBD}{\mbox{{\bf SBD}}}

\newcommand{\ecs}{\mbox{{\bf X}}}

\newcommand{\eg}{\mbox{{\bf g}}}

\newcommand{\paromega}{\partial \Omega}

\newcommand{\gau}{\Gamma_{u}}

\newcommand{\gaf}{\Gamma_{f}}

\newcommand{\sig}{{\bf \sigma}}

\newcommand{\gac}{\Gamma_{\mbox{{\bf c}}}}

\newcommand{\deu}{\dot{\eu}}

\newcommand{\dueu}{\underline{\deu}}

\newcommand{\dev}{\dot{\ev}}

\newcommand{\duev}{\underline{\dev}}

\newcommand{\weak}{\stackrel{w}{\approx}}

\newcommand{\mild}{\stackrel{m}{\approx}}

\newcommand{\lrightarrow}{\stackrel{L}{\rightarrow}}

\newcommand{\rrightarrow}{\stackrel{R}{\rightarrow}}

\newcommand{\strong}{\stackrel{s}{\approx}}

\newcommand{\weakdown}{\rightharpoondown}

\newcommand{\opg}{\stackrel{\mathfrak{g}}{\cdot}}

\newcommand{\opunu}{\stackrel{1}{\cdot}}
\newcommand{\opdoi}{\stackrel{2}{\cdot}}

\newcommand{\opn}{\stackrel{\mathfrak{n}}{\cdot}}
\newcommand{\opx}{\stackrel{x}{\cdot}}

\newcommand{\tr}{\ \mbox{tr}}

\newcommand{\Ad}{\ \mbox{Ad}}

\newcommand{\ad}{\ \mbox{ad}}

\renewcommand{\contentsname}{ }

\title{More than discrete or continuous: a bird's view}

\author{Marius Buliga \\
\\
Institute of Mathematics, Romanian Academy, 
P.O. BOX 1-764, \\
 RO 014700, Bucure\c sti, Romania\\
{\footnotesize Marius.Buliga@imar.ro}}

\date{This version:  19.11.2010}

\maketitle

\begin{abstract}
I try to give mathematical evidence to the following equivalence, which 
is based on ideas from Plato (Timaeus): reality emerges from a more primitive, non-geometrical, reality  in the same way as  the brain construct (understands, simulates, transforms, encodes or decodes) the image 
of reality, starting from intensive properties (like  a bunch of spiking signals sent by receptors in the retina), without any use of extensive (i.e. spatial or geometric)  properties. 
\end{abstract}

\vspace{.5cm}

\

\section{Main statement}

I shall try to give support to a statement based  on the following citation
 from Timaeus, by Plato  \cite{timaeus}. Further this will be called Plato hypothesis (or PH): 

\vspace{.5cm}

\paragraph{(PH)}"The sight in my opinion is the source of the greatest benefit to us, for had we never seen the stars, and the sun, and the heaven, none of the words which we have spoken about the universe would ever have been uttered. But now the sight of day and night, and the months and the revolutions of the years, have created number, and have given us a conception of time, and the power of enquiring about the nature of the universe; and from this source we have derived philosophy, than which no greater good ever was or will be given by the gods to mortal man." 

\vspace{.5cm}

\paragraph{Statement.}  Under Plato's hypothesis, it follows that: 
\begin{enumerate}
\item[(A)] reality emerges from a more primitive, non-geometrical, reality
\item[ ] in the same way as 
\item[(B)] the brain construct (understands, simulates, transforms, encodes or decodes) the image 
of reality, starting from intensive properties (like  a bunch of spiking signals sent by receptors in the retina), without any use of extensive (i.e. spatial or geometric)  properties. 
\end{enumerate}

The problem (B) is known in life sciences as the problem of "local sign" \cite{koen} \cite{koen2}. Indeed, any use of extensive 
properties would lead to the "homunculus fallacy" \cite{homunculus}.

These equivalent problems are difficult and wonderfully simple: 
\begin{enumerate}
\item[-] we don't know how to solve completely problem (A) (in physics) or problem (B) (in neuroscience),  
\item[-] but our ventral/dorsal streams and cerebellum do this all the time in about 150 ms. Moreover, 
any fly does it, as witnessed by its spectacular flight capacities  \cite{fly}. 
\end{enumerate}

The main argument in favor of (PH) comes from the fact that our knowledge in physics 
is rooted on the mathematical theories of calculus (Newton, Leibniz) and differential geometry (Riemann), which appeared  from a need to understand our physical world. Or, in section \ref{s3} is explained that maybe the (mathematical formalism of the)  differential  calculus and differential geometry could
  be based on the way the brain works (or biased by the way the brain works), specifically on the fact that we have binocular vision.  
  
The backbone of the argument is  that, as boolean logic is based on a primitive 
gate (like NAND), differential calculus and differrential geometry are also based on a primitive gate (a dilation gate), appearing naturally from the least 
sophisticated strategy of exploration which a binocular creature might have,  namely 
jumping randomly from a place to another, orienting herself by comparing what she 
sees with her two eyes.

\section{Orwellian, stalinesque or multiple drafts}

Is the reality discrete or continuous? This physics question,  related to problem (A), translates into the realm of problem (B) into the question: is the perceived reality discrete or continuous? Neuroscience has several intriguing answers to this. These answers can be then translated back into the language 
of problem (A). Let us see what we get in this way. 

Discrete sequences of events are perceived as continuous, as in the "beta illusion" \cite{betaillusion}, the more interesting "phi illlusion" \cite{phiillusion} or  the "cutaneous rabbit illusion" \cite{rabbit}. How is this possible? 
There are several theories explaining such illusions, nicely categorized 
by Daniel Dennet (also proposing his own theory of "multiple  drafts" \cite{dennet}). From the description given at \cite{denwiki}, such theories 
can be characterized as: 
\begin{enumerate}
\item[(a)] orwellian - "the subject comes to one conclusion, then goes back and changes that memory in light of subsequent events. This is akin to George Orwell's Nineteen Eighty-Four, where records of the past are routinely altered."
\item[(b)] stalinesque - the "events would be reconciled prior to entering the subject's consciousness, with the final result presented as fully resolved. This is akin to Joseph Stalin's show trials, where the verdict has been decided in advance and the trial is just a rote presentation."
\item[(c)] multiple drafts - "there are a variety of sensory inputs from a given event and also a variety of interpretations of these inputs". From  \cite{dennet2} 
[there is] "no central experiencer [who] confers a durable stamp of approval on any particular draft".
\end{enumerate}
Translated into the physics realm, this gives several interesting interpretations. 

\paragraph{(a)} Such a path has been pursued in physics, by Everett's 
 Many-Worlds Interpretation of Quantum Mechanics \cite{everett}. More precisely, concerning interpretations of the collapsing of the wave function which are  compatible with Everett theory, see Deutsch \cite{deut} and Stapp \cite{stapp}. 

\paragraph{(b)} In more general terms, not related especially to the problem of the discrete versus continuous nature of reality, we can see any theory based on extremality of an action like being of this type. However, probably due to my ignorance, I am not aware of physical theories supposing that a discrete reality conspires to give (to any observer) 
the appearance of being continuous. More precisely, such a theory would take as starting point a discrete reality  where discrete things happen, in the limit when the graininess goes to zero, like in a continuous reality. One big and  fundamental  difficulty would be then to give a reasonable mechanism of how is this possible. 

\paragraph{(c)} This could be seen in physics as if there is a "stack" of realities, one for 
every scale, say $reality(\varepsilon)$ at scale $\varepsilon > 0$, with the following property: if we take $reality(\varepsilon)$ and explore it at scale $\mu > 0$, then we 
get $reality(\varepsilon)(\mu)$, which is indiscernible from $reality(\varepsilon \mu)$.

\section{The metaphor of the binocular explorer} 
\label{s3}

I shall discuss about the simplest strategy to explore a space $X$ (say "a reality").  It may be useful to think about this in terms of distance and distance-preserving maps. 

\paragraph{Preliminaries.}Indeed, let us say that we are interested to know if a set $X$ is discrete or continuous. For this we send an explorer (call her Alice) to look around.   I shall suppose that we can put a distance on this set, that is a function 
$d: X \times X \rightarrow [0,+\infty)$ which satisfies the following requirement: 
we can make a map of any three points in $X$, in the sense that for any 
points $x,y,z \in X$ there is a bijective correspondence with a triple A,B,C in the plane  such that the sizes 
(lengths) of AB, BC, AC are equal respectively with $d(x,y)$, $d(y,z)$, $d(z,x)$. 
Basically, we accept that we can represent in the plane any three points from the 
space $X$. 

The explorer (call her Alice)  wants to make a map of 
a newly discovered land $X$ on the piece of paper $Y$ -- or by using the mound 
of clay $Y$,  or even by carefully recording in her collective mind the echos 
of her exploratory cries, if she were a mutant army of bats. 
"Understanding"  the space $X$ (with respect 
to the choice of the "gauge" function $d$) into the terms of the more familiar 
space $Y$ (itself endowed with a distance function $D$)  means  making a map of $X$ into $Y$ which is not deforming distances too much. Ideally, a perfect map has to be Lipschitz, that is the distances between points in $X$ are transformed by the map into distances between points in $Y$, with a precision independent of the scale: the map  $f:X \rightarrow Y$ is (bi)Lipschitz if there are positive numbers $c < C$ such that for any two points 
$x,y \in X$  
$$ c \, d(x,y) \, \leq \, D(f(x),f(y)) \, \leq \, C \, d(x,y)$$
This would be a very  good map, because it would tell how $X$ is at all scales.   Only that it is impossible to make such a map in practice. What we can do 
instead, is to sample the space $X$ (take a $\varepsilon$-dense subset of $X$ 
with respect to the distance $d$) and try to represent as good as possible this 
subspace in $Y$. Mathematically this is like asking for the map function $f$ to have the following property: there are supplementary positive constants $a, A$ such that 
for any two points 
$x,y \in X$  
$$ c \, d(x,y) - a  \, \leq \, D(f(x),f(y)) \, \leq \, C \, d(x,y) + A$$
The problem with these maps is that, only by using one  of them, we cannot decide 
if the space $X$ is discrete or continuous, because all small details are erased. 
There is hope though  to get an answer to our question  by using an "atlas" made by lots of such charts, such that for any $\varepsilon > 0$, there is a map 
with the constants $a/c, A/C$ smaller than $\varepsilon$.  Then, by a passage to the 
limit (in itself not at all clear), we should be able to probe the structure of 
$X$ at smaller and smaller scales. 

But is it possible that such an atlas exist, in the eventuality that the space $Y$ is not the same as space $X$? Mathematical results shows that not always!  This seems to be an issue that physicists ignore. 

\paragraph{The explorer's program.} Let us go into details of the following exploration program of the space $X$: the explorer Alice jumps randomly in the metric space $(X,d)$, making drafts of  maps at scale $\varepsilon$ and simultaneously orienting herself by using these draft maps. 

I shall explain each part of this exploration problem. 

\paragraph{Making maps at scale $\varepsilon$.} I shall suppose that Alice, while 
at $x \in X$, makes a map at scale $\varepsilon$ of a neighbourhood of 
$x$, called $\displaystyle V_{\varepsilon}(x)$, into another neighbourhood of 
$x$, called $\displaystyle U_{\varepsilon}(x)$. The map is a function: 
$$\delta^{x}_{\varepsilon} : U_{\varepsilon}(x) \rightarrow V_{\varepsilon}(x)$$ 
The "distance on the map" is just: 
\begin{equation}
d^{x}_{\varepsilon} (u,v) \, = \, \frac{1}{\varepsilon} d( \delta^{x}_{\varepsilon} u \, , \, \delta^{x}_{\varepsilon} v)
\label{drel}
\end{equation}
so the map $\displaystyle \delta^{x}_{\varepsilon}$ is indeed a rescaling map, with scale $\varepsilon$.

The choice of making maps of $X$ into $X$ is natural, because we may argue that the map of the reality is itself a part of the reality. 

We suppose of course that the map $\displaystyle \delta^{x}_{\varepsilon}$ is a bijection and we call the inverse of it by the name  $\displaystyle \delta^{x}_{\varepsilon^{-1}}$.  

We may suppose that the map of Alice, while being at point $x$, is a piece of paper 
laying down on the ground, centered such that the "$x$" on paper coincides with the place in $X$ marked by "$x$". In mathematical terms, I ask that 
\begin{equation}
\delta_{\varepsilon}^{x} x \, = \,  x
\label{comp0}
\end{equation}

\paragraph{Jumping randomly in the space $X$.} For random walks we need random  walk kernels, therefore we shall suppose that metric balls in $(X,d)$ have finite, non zero, (Hausdorff) measure. Let $\mu$ be  the associated  Hausdorff measure. The 
$\varepsilon$-random walk is then random  jumping from $x$ into the ball $B(x,\varepsilon)$, so the random walk kernel is  
\begin{equation}
m_{x}^{\varepsilon} \, = \,  \frac{1}{\mu(B(x,\varepsilon))} \,  \mu_{|_{B(x,\varepsilon)}} 
\label{drwk}
\end{equation}
but any Borel probability $\displaystyle m_{x}^{\varepsilon} \in \mathcal{P}(X)$ would be good. 

\paragraph{Compatibility considerations.} At his point I want to add some relations which give a more precise meaning to the scale $\varepsilon$.  

I shall first introduce a standard notation. If $\displaystyle f: X_{1} \rightarrow X_{2}$ is a Borel map and $\displaystyle \mu \in \mathcal{P}(X_{1})$ is a probability measure in 
the space $\displaystyle X_{1}$  then the push-forward of 
$\mu$ through $f$ is the probability measure $\displaystyle f \sharp \mu \in
\mathcal{P}(X_{2}$ defined by: for any $\displaystyle B\in\mathcal{B}(X_{2})$ 
$$\left(f \sharp \mu \right) (B) \, = \, \mu(f^{-1}(B))$$ 

For example,  notice that the random walk kernel is transported  into a random walk kernel by the inverse of the map $\displaystyle \delta^{x}_{\varepsilon^{-1}}$.  I shall impose  that for any (open) set $\displaystyle A \subset  U_{\varepsilon}(x)$ we have 
\begin{equation}
m_{x}^{\varepsilon}(\delta^{x}_{\varepsilon} (A)) \, = \, m_{x}(A) + \mathcal{O}(\varepsilon)
\label{comp1}
\end{equation}
where $\mathcal{O}(\varepsilon)$ is a function (independent on $x$ and $A$) going to zero as $\varepsilon$ is going to zero and $\displaystyle m_{x}$ is a random walk kernel of the map representation $\displaystyle (U_{\varepsilon}(x), d^{x}_{\varepsilon})$, more precisely this is a Borel probability. With the previously induced notation, relation (\ref{comp1}) appears as: 
$$\delta_{\varepsilon^{-1}}^{x} \, \sharp \, m_{x}^{\varepsilon} \, (A) \, = \, 
m_{x}(A) + \mathcal{O}(\varepsilon)$$

Another condition to be imposed is   that  $\displaystyle V_{\varepsilon}$ is approximately the ball $B(x,\varepsilon)$, in the following sense:
\begin{equation}
m_{x}(U_{\varepsilon}(x) \setminus \delta_{\varepsilon^{-1}}^{x} B(x,\varepsilon)) =  \mathcal{O}(\varepsilon)
\label{comp2}
\end{equation}

\paragraph{Multiple drafts reality.} We may see the  atlas that Alice draws as a  "multiple drafts" theory of $(X,d)$. Indeed, the reality at scale $\varepsilon$ 
(and centered at $x \in X$), according to Alice's exploration, can be seen as the triple: 
$$(U_{\varepsilon}(x), d^{x}_{\varepsilon}, \delta_{\varepsilon^{-1}}^{x} \, \sharp \, m_{x}^{\varepsilon})$$
Moreover, we may think that the data used to construct Alice's atlas are: the distance $d$ and for any $x \in X$ and $\varepsilon > 0$,  the functions $\displaystyle \delta_{\varepsilon}^{x}$ and the probability measures 
$\displaystyle m_{x}^{\varepsilon}$. 

For any $x \in X$ and $\varepsilon > 0$ these data may be transported to the 
"reality at scale $\varepsilon$, centered at $x \in X$". Indeed, transporting all 
the data by the function $\displaystyle \delta_{\varepsilon}^{x}$, we get 
$\displaystyle d^{x}_{\varepsilon}$ instead of $d$, and  for any 
$u \in X$ (sufficiently close to $x$ w.r.t. the distance $d$) and any  $\mu >0$,  
we define the relative map making functions: 
\begin{equation}
\delta^{x, u}_{\varepsilon, \mu} \, = \, \delta^{x}_{\varepsilon^{-1}} \, 
\delta_{\mu}^{\delta^{x}_{\varepsilon} u} \, \delta^{x}_{\varepsilon} 
\label{reldil}
\end{equation}
and the relative kernels of random walks: 
\begin{equation}
m_{x, u}^{\varepsilon, \mu} \, = \, \delta^{x}_{\varepsilon^{-1}} \, \sharp \, 
m_{u}^{\mu} 
\label{relrwk}
\end{equation}
With these data I define $reality(x,\varepsilon)$ as: 
\begin{equation}
reality(x,\varepsilon) \, = \, \left(x, \varepsilon, d^{x}_{\varepsilon} , (u,\mu) \mapsto \left( \delta^{x, u}_{\varepsilon, \mu} , m_{x, u}^{\varepsilon, \mu} \right) \right)
\label{real1}
\end{equation}
Remark that now we may repeat this construction and define $reality(x,\varepsilon)(y,\mu)$, and so on. A sufficient condition for having a "multiple drafts" like reality: 
$$reality(x, \varepsilon)(x, \mu) \, = \, reality(x, \varepsilon \mu)$$ 
is to suppose that for any $x \in X$ and any $\varepsilon , \mu > 0$ we have 
\begin{equation}
\delta_{\varepsilon}^{x} \,  \delta^{x}_{\mu} \, = \, \delta^{x}_{\varepsilon \mu}
\label{semigroup}
\end{equation}
For more discussion about this (but without considering random walk kernels) see 
section 2.3 \cite{buligaintro} and references therein about the notion of "metric profiles".

\paragraph{Binocular orientation.} Here I explain how Alice orients herself by using the draft maps from her atlas. Remember that Alice jumps from $x$ to $\displaystyle \delta_{\varepsilon}^{x} u$ (the point $u$ on her map centered at $x$). Once arrived at 
$\displaystyle \delta_{\varepsilon}^{x} u$, Alice draws another map and then she uses a binocular approach to understand what she sees from her new location. She compares the maps 
simultaneously. Each map is like a telescope, a microscope, or an eye. Alice has two eyes. With one she sees $reality(x,\varepsilon)$ and, with the other, $\displaystyle 
reality(\delta_{\varepsilon}^{x}u , \varepsilon)$. These two "retinal images" are compatible in a very specific sense, which gives to Alice 
the sense of her movements. Namely Alice uses the following mathematical fact, which was explained before several times, in the language of: "dilatation structures" \cite{buligadil1} \cite{buligadil2}, "metric spaces with dilations" \cite{buligaintro} (but not using groupoids), 
"emergent algebras" \cite{buligairq} (showing that the distance is not necessary 
for having the result, and also not using groupoids), \cite{buligagr} (using 
normed groupoids). 

\begin{theorem}
There is a groupoid (a small category with invertible arrows) $Tr(X,\varepsilon)$ which has as objects $reality(x,\varepsilon)$, defined at (\ref{real1}), for all 
$x \in X$, and  as arrows the functions: 
\begin{equation}
\Sigma_{\varepsilon}^{x}(u, \cdot) : \, reality(\delta_{\varepsilon}^{x}u , \varepsilon) \rightarrow \, reality(x,\varepsilon) 
\label{real2}
\end{equation}
defined for any $\displaystyle u, v \in U_{\varepsilon}(x)$ by: 
\begin{equation}
\Sigma_{\varepsilon}^{x}(u,v) \, = \, \delta^{x}_{\varepsilon^{-1}} \, \delta^{\delta^{x}_{\varepsilon} u}_{\varepsilon}  v
\label{sig1}
\end{equation}
 Moreover, all arrows are 
isomorphisms (they are isometries and transport in the right way the structures defined at (\ref{reldil}) and (\ref{relrwk}), from the source of the arrrow to the 
target of it). 
\label{t1}
\end{theorem}

Therefore Alice, by comparing the maps she has, orients herself by using the 
"approximate translation by $u$" 
$\displaystyle \Sigma_{\varepsilon}^{x}(u, \cdot)$ 
(\ref{sig1}), which appears as an isomorphism between the two realities, cf. (\ref{real2}). This function deserves the name "approximate translation", as proved elsewhere \cite{buligadil1} \cite{buligadil2} \cite{buligairq}, because of 
the following reason: if we suppose that $reality(x,\varepsilon)$ converges 
as the scale $\varepsilon$ goes to $0$, then $\displaystyle \Sigma_{\varepsilon}^{x}(u, \cdot)$ should converge to the translation by $u$, as seen  in the tangent space as $x$. Namely we have the following mathematical definition and theorem  (for the proof see the cited references; the only new thing concerns the convergence of the random walk kernel, but this is an easy consequence of the compatibility conditions (\ref{comp1}) (\ref{comp2})). 

\begin{definition}
A metric space with dilations  and random walk (or dilatation structure with a random walk) is a structure $(X,d,\delta,m)$ such that:  
\begin{enumerate}
\item[(a)] $(X,d,\delta)$ is a normed uniform idempotent right quasigroup  (or  dilatation structure), 
cf. \cite{buligabraided} Definition 7.1. 
\item[(b)] $m$ is a random walk kernel, that is for every $\varepsilon > 0$ 
we have a measurable function $\displaystyle x \in X \, \mapsto \, m_{x}^{\varepsilon} \in \mathcal{P}(X)$ (a transversal function on the pair groupoid 
$X \times X$, in the sense of \cite{connes} p. 35) which satisfies the compatibility conditions (\ref{comp1}) (\ref{comp2}). 
\end{enumerate}
\label{de1}
\end{definition}

\begin{theorem}
Let $(X,d,\delta,m)$ be a   dilatation structure with a random walk. Then for any 
$x \in X$, as the scale $\varepsilon$ goes to $0$ the structure $reality(x,\varepsilon)$ converges to $reality(x,0)$ defined by: 
\begin{equation}
reality(x,0) \, = \, \left(x, 0, d^{x} , (u,\mu) \mapsto \left( \delta^{x, u}_{ \mu} , m_{x, u}^{ \mu} \right) \right)
\label{realtang}
\end{equation}
in the sense that $\displaystyle d^{x}_{\varepsilon}$ converges uniformly 
 to $\displaystyle d^{x}$, dilations $\displaystyle \delta^{x, u}_{\varepsilon,  \mu}$ converge uniformly to dilations $\displaystyle \delta^{x, u}_{ \mu}$, 
 probabilities  $\displaystyle m_{x, u}^{\varepsilon, \mu}$  converge simply  
 to probabilities $\displaystyle  m_{x, u}^{ \mu}$. Finally $\displaystyle 
 \Sigma_{\varepsilon}^{x}(\cdot , \cdot)$ converges uniformly to a conical group 
 operation. 
 \label{t2}
 \end{theorem}
 
Conical groups are a (noncommutative and vast) generalization of real vector 
spaces. See for this \cite{buligadil2}. 

\paragraph{Alice knows the differential geometry of $X$.} All is hidden in dilations (or maps made by Alice) $\displaystyle \delta_{\varepsilon}^{x}$. They encode the approximate (theorem \ref{t1}) and 
exact (theorem \ref{t2}) differential AND algebraic structure of the tangent 
bundle of $X$.  Moreover, if Alice feels the need to differentiate, then again 
she uses the maps she has (dilations) to rescale the function she wants to differentiate (see \cite{buligaspace} for a graphical,  but not geometrical,   interpretation of this formalism by using braids diagrams).

\section{Discussion}

\paragraph{The metaphor of the universal front end.}  Maybe biology uses at a different scale an embodiment of a fundamental mechanism of the nature.  
Indeed, I suggest  that  the physical space can be understood, at some very fundamental level, as the input of a "universal front end", and physical 
observers are "universal front ends".

Then, as suggested in \cite{buligaspace},  physical observers are like universal front ends looking at the same (but otherwise unknown) space. This is compatible with the   main statement of the paper, which then  suggests that there might be an analogy between the problem of "local sign" in neuroscience and  the problem of understanding the properties of the physical space as emerging from some non-geometrical, more  fundamental structure, like a net, a foam, a graph... 

I see no obstacle for giving the same  treatment to abstract mathematics. We may imagine that there is an abstract mathematical "front end" which, if fed with the definition of  a "space",  then spews out a "data structure" which is used for "past processing", that is for mathematical reasoning in that space. (In fact, when we say "let $M$ be a manifold", for example, we don't "have" that manifold,  only some properties of it, together with some very general abstract nonsense concerning  "legal" manipulations in the universe of "manifolds". All these can be liken with the image that we get 
 past the "front end" , in the sense that, like a real perceived image, we see it all, but we are incapable 
 of really enumerating and naming all that we see.)

\paragraph{Alice communicates with Bob.} Now Alice wants to send Bob her knowledge. 
Bob lives in the metric space $(Y,D)$, which was well explored previously, therefore 
himself knows well the differential geometry and calculus in the space $Y$. 

Based on Alice's informations, Bob hopes to construct a Lipschitz function from an open set in $X$ with image an open set in $Y$. If he succeeds, then he will know all that  Alice knows (by transport of all relevant structure using his Lipschitz function). But he might fail, because of a strong mathematical result (unknown apparently to physicists, apparently  aware only about the dichotomy discrete/continuous). Indeed, as a consequence of Rademacher theorem (see Pansu \cite{pansu}, for a 
Rademacher theorem for Lipschitz function between Carnot groups, a type of conical groups), if such a function would exist then it would be differentiable almost everywhere. Or, the differential (at a point) is a morphism of conical groups which 
commutes with Alice's and Bob's maps (their respective dilations). Or, it might happen that there is no non-trivial such morphism and we arrrive at a contradiction! For example, it is  known that there is no Lipschitz embedding 
of an  open set in Heisenberg group with its Carnot-Carath\'eodory distance into 
the Hilbert space. So, if Alice explored the Heisenberg group and if Bob expects to understand in Hilbert space, or even in $L^1$, what Alice saw, then it fails (in the sense explained previously). See for this \cite{cheeger} and references therein.  

This is not a matter of topology (continuous or discrete, or 
complex topological features at all scales), but a matter having to do with the 
construction of the "realities" (in the sense of relation (\ref{real1}), having  all to do with differential calculus and differential geometry in the most fundamental sense.

 Therefore, another strategy of Alice would be to communicate to Bob not all details of her map, but the relevant algebraic identities that her maps (dilations) satisfy. Then Bob may try to simulate what Alice saw. This is  a path first suggested in  \cite{buligaspace}.

\section{Explanation of the title, according to Plato, Timaeus}

"But the race of birds was created out of innocent light-minded men, who, although their minds were directed toward heaven, imagined, in their simplicity, that the clearest demonstration of the things above was to be obtained by sight; these were remodelled and transformed into birds, and they grew feathers instead of hair." (source 
\cite{timaeus})

\section{Appendix 1: Plato, Timaeus}

The ancient greeks understanding of vision can be inferred from the dialogue Timaeus, by Plato. Here are some excerpts, where this is explained in a concise manner (source  \cite{timaeus}).

\vspace{.5cm}

Timaeus, 12:  "First, then, the gods, imitating the spherical shape of the universe, enclosed the two divine courses in a spherical body, that, namely, which we now term the head, being the most divine part of us and the lord of all that is in us..."  

"And of the organs they first contrived the eyes to give light, and the principle according to which they were inserted was as follows: So much of fire as would not burn, but gave a gentle light, they formed into a substance akin to the light of every-day life; and the pure fire which is within us and related thereto they made to flow through the eyes in a stream smooth and dense, compressing the whole eye, and especially the centre part, so that it kept out everything of a coarser nature, and allowed to pass only this pure element."

"And the whole stream of vision, being similarly affected in virtue of similarity, diffuses the motions of what it touches or what touches it over the whole body, until they reach the soul, causing that perception which we call sight."

"And now there is no longer any difficulty in understanding the creation of images in mirrors and all smooth and bright surfaces. For from the communion of the internal and external fires, and again from the union of them and their numerous transformations when they meet in the mirror, all these appearances of necessity arise, when the fire from the face coalesces with the fire from the eye on the bright and smooth surface. And right appears left and left right, because the visual rays come into contact with the rays emitted by the object in a manner contrary to the usual mode of meeting..."

\section{Appendix 2: Computations in the front end visual system as a paradigm}

Koenderink, Kappers and van Doorn \cite{koen2} define the "front end visual system" as "the interface between the light field and those parts of the brain nearest to the transduction stage". The authors propose the following 
characterization of the front end visual system, section 1 \cite{koen2}: 
\begin{enumerate}
\item[1.] "the front end is a "machine" in the sense of a syntactical transformer; 
\item[2.] there is no semantics. The front end processes structure;
\item[3.] the front end is a deterministic machine; ... all output depends causally on the (total) input from the immediate past." 
\item[4.] "What is not explicitly encoded by the front end is irretrievably lost. Thus the front end should be 
universal (undedicated) and yet should provide explicit data structures."
\end{enumerate}

 Koenderink et al.  propose (some embodiment of differential) geometry as a model  for this activity in the brain, \cite{koen}, p. 126: 
\vspace{.5cm}

"The brain can organize {\em itself} through information obtained via interactions with the physical world into 
an embodiment of geometry, it becomes a veritable {\em geometry engine} [...] 
many of the better-known brain processes are most easily understood in terms of differential geometrical calculations running on massively parallel processor arrays whose nodes can be understood quite directly in terms of multilinear operators (vectors, tensors, etc). In this view brain processes in fact are space."

\vspace{.5cm}

 The front end functions in a  massively parallel manner, thus it has to do it by working with local representations. I cite from the last paragraph of section 1, just before the beginning of sections 1.1 \cite{koen2}.

\vspace{.5cm}

"In a local representation one can do without extensive (that is spatial, or geometrical) properties and represent everything in terms of intensive properties. This obviates the need for explicit geometrical expertise. The local representation of geometry is the typical tool of differential geometry. ... The columnar organization of representation in primate visual cortex suggests exactly such a structure."


\begin{thebibliography}{99}

\bibitem{timaeus} The Internet Classics Archive, Timaeus, by Plato, translated by Benjamin Jowett, available at:  \url{http://classics.mit.edu//Plato/timaeus.html}

\bibitem{betaillusion} \url{http://en.wikipedia.org/wiki/Beta_movement}

\bibitem{rabbit} \url{http://en.wikipedia.org/wiki/Cutaneous_rabbit_illusion}

\bibitem{phiillusion} \url{http://en.wikipedia.org/wiki/Phi_phenomenon}

\bibitem{homunculus} \url{http://en.wikipedia.org/wiki/Homunculus_argument}


\bibitem{buligadil1} M. Buliga, Dilatation structures I. Fundamentals, {\it 
J. Gen. Lie Theory Appl.},  {\bf 1} (2007),  2, 65-95. 

\bibitem{buligaintro} M. Buliga, Introduction to metric spaces with dilations, (2010),   \url{http://arxiv.org/abs/1007.2362}

\bibitem{buligadil2} M. Buliga, Infinitesimal affine geometry of metric spaces 
endowed with a dilatation structure ,  {\it Houston Journal 
of Math.} 36, 1 (2010), 91-136,  \url{http://arxiv.org/abs/0804.0135}


\bibitem{buligairq} M. Buliga, Emergent algebras as generalizations of
differentiable algebras, with applications, (2009), \url{http://arxiv.org/abs/0907.1520}

\bibitem{buligaspace} M. Buliga, What is a space? Computations in emergent algebras and the front end visual system, (2010) \\ 
 \url{http://arxiv.org/abs/1009.5028}

\bibitem{buligabraided} M. Buliga, Braided spaces with dilations and sub-riemannian symmetric spaces, (2010),  \url{http://arxiv.org/abs/1005.5031} 

\bibitem{buligagr} M. Buliga, Deformations of normed groupoids and differential calculus. 
First part, (2009), \url{http://arxiv.org/abs/0911.1300}

\bibitem{cheeger} J. Cheeger, B. Kleiner, Metric differentiation, monotonicity and maps to $L^1$, \url{http://arxiv.org/abs/0907.3295}

\bibitem{connes} A. Connes, Sur la theorie non commutative de l'integration, in: Alg\'ebres d'Op\'erateurs, S\'eminaire sur les Alg\'ebres d'Op\'erateurs, Les Plans-sur-Bex, Suisse, 13-18 mars 1978, Lecture Notes in Mathematics 725, ed. by 
A. Dold and B. Eckmann, Springer-Verlag 1079, p. 19-143


\bibitem{dennet} D.C. Dennett, (1991) Consciousness Explained. Boston: Little Brown

\bibitem{dennet2} D.C. Dennett, M. Kinsbourne. "Time and the Observer: the Where and When of Consciousness in the Brain". Behavioral and Brain Sciences (15): 183–247, (1992)

\bibitem{denwiki} \url{http://en.wikipedia.org/wiki/Multiple_Drafts_Model}


\bibitem{deut} D. Deutsch, Quantum theory as a universal physical theory, 
{\it Int. J. Theor. Phys.} {\bf 24}, 1-41 (1985)

\bibitem{everett} H. Everett, Theory of the Universal Wavefunction, Thesis, Princeton University, (1956, 1973), pp 1-140



\bibitem{fly} N. Franceschini, J.M. Pichon, C. Blanes, From insect vision to 
robot vision, {\it Phil. Trans.: Biological Sciences},  {\bf 337},   1281 (1992), 
Natural and Artificial Low-Level Seeing Systems, 283-294




\bibitem{koen} J. Koenderink, The brain a geometry engine, {\it Psychol. Res.} {\bf 52} (1990), 122-127

\bibitem{koen2} J.. Koenderink, A. Kappers, A. van Doorn, Local Operations :The Embodiment of Geometry. Basic Research Series, (1992), 1-23



\bibitem{pansu} P. Pansu, M\'etriques de Carnot-Carath\'eodory et
quasi-isom\'etries des espaces sym\'etriques de rang un, Ann. of Math., (2) 
{\bf 129}, (1989), 1-60



\bibitem{stapp} H.P. Stapp, The basis problem in many-worlds theories, {\it 
Can. J. Phys.} {\bf 80}(9), 1043–1052 (2002) 




\end{thebibliography}
\end{document}